\documentclass[12pt]{article}
\usepackage{amssymb}
\usepackage{hyperref}
\usepackage{latexsym}
\usepackage{color}
\usepackage{amsmath}
\usepackage{graphicx}
\usepackage{amsfonts}
\usepackage{amsthm}
\usepackage{graphicx}
\usepackage{mathrsfs}
\newtheorem{theorem}{Theorem}
\newtheorem{lemma}{Lemma}

\newtheorem{remark}{Remark}

\everymath{\displaystyle}

\begin{document}

\newcommand{\al}{\alpha}
\newcommand{\bt}{\beta}
\newcommand{\ep}{\varepsilon}
\newcommand{\dl}{\delta}
\newcommand{\ti}{\tilde}
\newcommand{\be}{\begin{equation}}
\newcommand{\ee}{\end{equation}}
\newcommand{\goto}{\rightarrow}
\newcommand{\half}{\frac{1}{2}}
\newcommand{\ld}{\lambda}
\newcommand{\p}{\partial}
\newcommand{\ph}{\varphi}
\newcommand{\ov}{\overline}
\newcommand{\pl}{\parallel}
\newcommand{\br}{\breve}
\newcommand{\var}{\varphi}

\newcommand\Ker{{\rm Ker}}
\newcommand\intl{\int\limits}
\newcommand\suml{\sum\limits}
\newcommand\maxl{\max\limits}
\newcommand\minl{\min\limits}
\newcommand\supl{\sup\limits}
\newcommand\infl{\inf\limits}
\newcommand\liml{\lim\limits}
\newcommand\ch{{\rm ch}}
\newcommand\tg{{\rm tg}}
\newcommand\rank{{\rm rank}}
\newcommand\const{{\rm const}}
\newcommand\Ga{\Gamma}
\renewcommand{\Im}{{\rm Im}}
\makeatletter \@addtoreset{equation}{section}
\renewcommand{\theequation}{\thesection.\arabic{equation}}
\makeatother

\font\bfsl=cmssi17 \font\sectiontt=cmtt12 scaled\magstep1
\font\stt=cmtt10 scaled 850

\renewcommand{\baselinestretch}{1}
\renewcommand{\theequation}{\arabic{section}.\arabic{equation}}
\renewcommand{\thetheorem}{\arabic{section}.\arabic{theorem}}
\renewcommand{\thelemma}{\arabic{section}.\arabic{lemma}}
\renewcommand{\theproposition}{\arabic{section}.\arabic{proposition}}
\renewcommand{\thedefinition}{\arabic{section}.\arabic{definition}}
\renewcommand{\thecorollary}{\arabic{section}.\arabic{corollary}}
\renewcommand{\theremark}{\arabic{section}.\arabic{remark}}
\pagestyle{empty}

\pagestyle{myheadings}

\title{Lower Bound For The Ratios Of Eigenvalues
Of Schr\"{o}dinger Equations With Nonpositive Single-Barrier
Potentials}

\author{Jamel Ben Amara \thanks{Faculty of Sciences of Tunis,
University of Tunis El Manar, Mathematical Engineering Laboratory,
Polytechnic School of Tunisia,
jamel.benamara@fsb.rnu.tn}~~~~~~~~~~Jihed Hedhly
\thanks{Faculty of Sciences  Bizerte, University of Carthage,
  Mathematical Engineering Laboratory,
Polytechnic School of Tunisia, hjihed@gmail.com.}}

\date{}

 \maketitle

\begin{abstract} Horv\'{a}th and Kiss [Proc. Amer. Math. Soc., 2005] proved the upper
bound estimate  $\frac{\lambda _{n}}{\lambda _{m}}\leq
\frac{n^{2}}{m^{2}}$ $ (n>m\geq 1) $ for Dirichlet eigenvalue ratios
of the Schr\"{o}dinger problem $-y''+q(x)y=\lambda y$ with
nonnegative and single-well potential $q$. In this paper, we prove
that if $q(x)$ is a nonpositive, continuous and single-barrier
potential, then $\frac{\lambda_{n}}{\lambda_{m}}\geq
\frac{n^{2}}{m^{2}}$ for $\lambda_n>\lambda_m \geq -2q^*$, where
$q^{\ast}=\min\{q(0), q(1)\}$. In particular, if $q(x)$ satisfies
the additional condition $\mid q^{\ast } \mid\leq
\frac{\pi^{2}}{3}$, then $\lambda _{1}>0$ and
$\frac{\lambda _{n}}{\lambda _{m}}\geq \frac{n^{2}%
}{m^{2}}$ for $n>m\geq 1.$ For this result, we develop a new
approach to study the monotonicity of the modified Pr\"{u}fer angle
function. 
\end{abstract}
~~~~~~~~~~
\maketitle
\section{ Introduction}
Consider the  One-dimensional Schr\"{o}dinger equation acting on $[0,1]$%
\begin{equation} \label{1.1}
-y^{\prime \prime }+q(x)y=\lambda y,
\end{equation}%
with Dirichlet boundary conditions
\begin{equation}
y(0)=y(1)=0.  \label{1.2}
\end{equation} Here $q$ is a nonpositive continuous and single-barrier
potential in $[0,1].$

Note that, the string equation $-y''=\lambda \rho(x) y,$ with $\rho$
is twice differentiable and $\rho^{-\frac{1}{4}}$ is convex, can be
transformed into \eqref{1.1} (where $q\leq0$) by the Liouville
substitution (see \cite {8}).

It is known  (see \cite{8})  that the spectrum of Problem
\eqref{1.1} - \eqref{1.2} consists of a growing sequence of
infinitely point $\lambda _{1}<\lambda _{2}<..........<\lambda
_{n}...\infty$.

The issues of optimal estimates for the eigenvalue ratios
$\frac{\lambda_n}{\lambda_m}$ have attracted a lot of attention
(cf.\cite{1,2',JJ,H2016,MJH,5})
and references therein. In $1989$, Ashbaugh and Benguria $\cite{1}$
proved the optimal upper bound $\frac{\lambda_{n}}{\lambda_{1}}\leq
n^{2},$ for nonnegative potentials. 
 In $1996$, Huang and Law
$\cite{3}$ extended the results in $\cite{1}$ to more general
boundary conditions. Chen et al. in $\cite{4}$
proved the lower bound $\frac{\lambda _{n}}{\lambda _{m}}\geq ( \lfloor \frac{n}{m}%
\rfloor ) ^{2},$ for nonpositive potentials, where $\lfloor s
\rfloor$ denotes the largest integer less than or equal to $s$. In
$1998$, Law and Huang in \cite{CK} proved that the eigenvalues of
the regular Sturm-Liouville equation $-(p(x)y^{\prime})^ {\prime
}+q(x)y=\lambda \rho (x)y$ (with Dirichlet boundary conditions)
satisfy the lower bound
$$\frac{\lambda_{n}}{\lambda_{m}}\geq\frac{1}{1+\xi}
\Big(\frac{n+1}{m+1}\Big)^2\frac{k}{K},~~n>m\geq0,$$ for $p\in
C^1[0,1],$ $q,\rho\in C[0,1]$, $q\geq0$ and $0<k\leq p(x)\rho(x)\leq
K,$ with $\xi=\frac{K\max{\{pq\}}}{k(m+1)^2\sigma^2\pi^2}$ and
$\sigma=\Big(\int_{0}^{1}\frac{1}{p(s)}ds\Big)^{-1}$. In $2005$,
Horv\'{a}th and Kiss $\cite{5}$ showed that if $q(x)$ is a
nonnegative single-well potential, then \begin{eqnarray}\label{33}
\frac{\lambda _{n}}{\lambda _{m}}\leq \frac{n^{2}}{m^{2}},~~ n>
m\geq1.\end{eqnarray} Their approach is mainly based on the
monotonicity of the Pr\"{u}fer angle as function in $\lambda>0$ (see
$\cite[~Theorem ~2.2.]{5}$). At the end of their paper
$\cite[~Remark~5.1.]{5}$, they gave an example of a single-barrier
potential which shows that the Pr\"{u}fer angle is not a monotonous
function. Recently, the authors \cite{JJ3} proved \eqref{33} for one
class of nonnegative differentiable and single-barrier potentials.

A function in $[0, 1]$ is called single-barrier (single-well) if
there is a point $x_0$ in $[0, 1]$, such that it is monotone
increasing (decreasing) in $[0,x_0]$ and monotone decreasing
(increasing) in $[x_0,1]$ (see in \cite{2}).

In the present paper, we prove that if $q(x)$ is a nonpositive,
continuous and single-barrier potential, then
 $\frac{\lambda _{n}}{\lambda _{m}}\geq \frac{n^{2}%
}{m^{2}}$ for $\lambda _{n}>\lambda _{m}\geq -2q^{\ast },$ where
$q^{\ast}=\min\{q(0), q(1)\}$. In particular, if $q(x)$ satisfies
the additional condition $|q^{\ast}|\leq \frac{\pi^{2}}{3}$, then
$\lambda _{1}>0$ and $\frac{\lambda _{n}}{\lambda _{m}}\geq
\frac{n^{2}}{m^{2}},$ for $n>m\geq 1$. For this result, we prove
that the associated Pr\"{u}fer angle is a decreasing function in
$\lambda>0.$ Note that (see Remark \ref{rem1}) our approach used in
this paper can be applied to the case of nonnegative single-well
potentials studied in $\cite{5}$.
\section{Preliminaries And The  Main Statements} \label{c}
Denote by $y( x,z ) $ the unique solution of the initial value
problem \begin{eqnarray} \label{H.S} \left\{
\begin{array}{ll}-y^{\prime \prime }+q(x)y=z^{2}y,\ \ \ x\in [ 0,1] ,\ \ \ z>0,\\
y(0)=0, \ \ \  y^{\prime }(0)=1.
 \end{array}
 \right.
\end{eqnarray}
We shall apply to System \eqref{H.S}, the modified Pr\"{u}fer
substitution as introduced in $\cite{1}$. Let
\begin{eqnarray}\label{2.2}
&&y(x,z) =r( x,z) \sin \varphi ( x,z), \cr &&y^{\prime }(
x,z)=zr(x,z) \cos \varphi ( x,z), \cr &&\varphi (0,z) =0,
\end{eqnarray}
where $r(x,z)>0$, and then let $\theta(x,z) =\frac{\varphi (x,z)}{z}.$\\
Using Equation $\eqref{1.1} $ together with \eqref{2.2}, one finds
the following differential equations for $r( x,z)$ and
$\varphi(x,z)$:
\begin{equation} \label{2.6}
\varphi ^{\prime }=z-\frac{q}{z}\sin ^{2}\varphi ,
\end{equation}%
\begin{equation}\label{2.7}
\frac{r^{\prime }}{r}=\frac{q}{z}\sin \varphi \cos \varphi .
\end{equation}%
It is obvious that $z^{2}$ is an eigenvalue iff $\varphi(\pi ,z) $
is a multiple of $\pi $. Denote by $z_{n}$ the square root of the
eigenvalue $\lambda _{n}$ of \eqref{1.1}-\eqref{1.2}. Since $q(x)$
is nonpositive then by \eqref{2.6}, $\varphi ^{\prime}>0$ for $z>0.$
In this case $\varphi ^{-1}$ exists and $\varphi
^{-1}(k\pi+\frac{\pi}{2})$, $\varphi ^{-1}((k+1)\pi)$ ($k\in
\mathbb{N}$) are the zeros of $y'$ and $y$ in $(0,x_{0}]$,
respectively. It is known (e.g., see \cite[~chap.1]{5}) that these
zeros are decreasing as $z$ increases. We denote by prime $({resp.~
dot})$ the derivative with respect to $x$ $({resp. ~z})$.  We now
enunciate the main results of this paper.

\begin{theorem}\label{the1}
Let $q(x)\leq 0$ be a continuous monotone increasing potential in
$[0,x_{0}].$ Then $\dot{\theta}( x_{0},z) \leq 0$ for\ $z\geq
\sqrt{-2q( 0) }.$ If there is a $z\geq \sqrt{-2q( 0) }$ with $\dot{
\theta}( x_{0},z) =0,$ then $q\equiv 0$ in $[0,x_{0}]$.
\end{theorem}
The proof will be given in Section \ref{a}.
\begin{theorem}\label{the2}
For the Schr\"{o}dinger problem $\eqref{1.1}-\eqref{1.2} ,$ if
$q(x)\leq 0$ is a continuous and single-barrier potential, then
\begin{equation}
\frac{\lambda _{n}}{\lambda _{m}}\geq \frac{n^{2}}{m^{2}},~for~
\lambda _{n}>\lambda _{m}\geq -2q^{\ast },
\end{equation}
where $q^{\ast}=\min\{q(0), q(1)\}$. In particular, if $q(x)$
satisfies the additional condition $|q^{\ast}|\leq
\frac{\pi^{2}}{3}$, then $\lambda _{1}>0$ and
$$\frac{\lambda _{n}}{\lambda _{m}}\geq \frac{n^{2}
}{m^{2}},~ for ~n>m\geq 1.$$ If for two different $m$ and $n$ the
equality holds, then $q\equiv 0 $ in $[0,1].$
\end{theorem}
The proof of Theorem \ref{the2} will be given in Section \ref{b}.

\section{The Proofs Of Theorems \ref{the1} }  \label{a}
For the proof of Theorem \ref{the1} we need the following results.

\begin{lemma}\label{cor1}(Corollary 3.3 in \cite{5})
\begin{equation}\label{2.9}
\dot{\theta}( x,z) =\frac{2}{z^{2}r^{2}(x)}\int_{0}^{x}r^{2}(t)%
\frac{q(t)}{z}\Big( \sin ^{2}\varphi -\varphi \sin \varphi \cos
\varphi \Big) dt.
\end{equation}
\end{lemma}
\begin{lemma} \label{lem3}  (Lemma 3.4 in \cite{5})
If $\vert \varphi \vert \in ] 0,\frac{%
\pi }{2}[ ,$ then $\sin ^{2}\varphi -\varphi \sin \varphi \cos
\varphi
>0.$
\end{lemma}
\begin{lemma}\label{lem4}
Let $k\geq 0$ be an integer and $0\leq C\leq \frac{\pi }{2
},$ $0\leq D\leq \pi,$ then%
\begin{equation} \label{A.A}
\int_{0}^{\varphi ^{-1}( C) }r^{2}(x)\frac{q(x)}{z}\Big( \sin
^{2}\varphi -\varphi \sin \varphi \cos \varphi \Big) dx \leq 0,
\end{equation}%
\begin{eqnarray} \label{A.C}
&&\int_{\varphi ^{-1}( k\pi +\frac{\pi }{2}) }^{\varphi ^{-1}( k\pi
+\frac{\pi }{2}+D) }r^{2}(x)\frac{q(x)}{z}\Big(\sin ^{2}\varphi
-\varphi \sin \varphi
\cos \varphi\Big) dx \cr &&\leq -( k+1) \pi \int_{\varphi ^{-1}( k\pi +\frac{\pi }{2}%
) }^{\varphi ^{-1}( k\pi +\frac{\pi }{2}+D) }r^{2}(x)\frac{%
q(x)}{z}\sin \varphi \cos \varphi dx,
\end{eqnarray}%
and equality holds iff $q\equiv0$ in the corresponding interval.

\end{lemma}
\begin{proof}
The first inequality $\eqref{A.A}$ follows from Lemma \ref{lem3}.
Following the proof of Corollary 3.5 in \cite{5}, and using the fact
that $|\varphi -( k+1) \pi |$ $\in]0,\frac{\pi }{2}[$, we obtain
\begin{eqnarray*}
&&\int_{\varphi ^{-1}( k\pi +\frac{\pi }{2}) }^{\varphi ^{-1}( k\pi
+\frac{\pi }{2}+D) }r^{2}( x) \frac{q( x) }{z}\Big( \sin ^{2}\varphi
-\varphi \sin \varphi \cos
\varphi \Big) dx \\
&=&\int_{\varphi ^{-1}( k\pi +\frac{\pi }{2}) }^{\varphi
^{-1}( k\pi +\frac{\pi }{2}+D) }r^{2}( x) \frac{%
q( x) }{z}\Big( \sin ^{2}\varphi -[ \varphi -(
k+1) \pi ] \sin \varphi \cos \varphi \Big) dx \\
&&-( k+1) \pi \int_{\varphi ^{-1}( k\pi +\frac{\pi }{2}%
) }^{\varphi ^{-1}( k\pi +\frac{\pi }{2}+D) }r^{2}(x)\frac{%
q(x)}{z}\sin \varphi \cos \varphi dx \\
&\leq &-( k+1) \pi \int_{\varphi ^{-1}( k\pi +\frac{\pi }{2}%
) }^{\varphi ^{-1}( k\pi +\frac{\pi }{2}+D) }r^{2}(x)\frac{%
q(x)}{z}\sin \varphi \cos \varphi dx.
\end{eqnarray*}
If equality holds, then $$\int_{\varphi ^{-1}( k\pi +\frac{\pi }{2})
}^{\varphi ^{-1}( k\pi +\frac{\pi }{2}+D) }r^{2}( x) \frac{q(x)
}{z}\Big(\sin ^{2}\varphi -[ \varphi -( k+1) \pi ] \sin \varphi \cos
\varphi \Big) dx =0,$$ whence $q\equiv0.$
\end{proof}
\begin{lemma}\label{lem5}
If $q(x)\leq 0,$ increasing in $[ 0,x_{0}] $ and $z^{2}\geq
-2q( 0) $ , then for each integer $n\geq 0,$ the function $%
x\mapsto  q^{2n+1}(x)+\frac{q^{2n+2}(x)}{z^{2}}$ is increasing in
$[0,x_{0}]$.
\end{lemma}
\begin{proof}
Let $n=0,$ and $x_{1},x_{2}\in [ 0,x_{0}] $ with $x_{1}<x_{2}. $
Then
\begin{eqnarray*}
&&q( x_{2}) +\frac{q^{2}( x_{2}) }{z^{2}}-q( x_{1}) -\frac{q^{2}(
x_{1}) }{z^{2}} =\Big( q( x_{2}) -q( x_{1})\Big) +\frac{1}{z^{2}}
\Big( q^{2}( x_{2}) -q^{2}( x_{1})\Big) \\
&=&\Big(q( x_{2}) -q( x_{1})\Big) \Big( 1+\frac{1}{ z^{2}}( q(
x_{2}) +q( x_{1}) )\Big) \geq \Big(q( x_{2}) -q( x_{1})\Big) \Big(
1+\frac{2q( x_{1}) }{z^{2}}\Big) \\&&\geq \Big(q( x_{2}) -q(
x_{1})\Big) \Big( 1+\frac{2q(0) }{z^{2}}\Big)\geq0 .
\end{eqnarray*}%
Suppose that $x\mapsto q^{2n+1}+\frac{q^{2n+2}}{z^{2}}$ is
increasing in
$x\in[ 0,x_{0}] $ for $n\geq 1.$ Then%
\begin{eqnarray*}
&&q^{2n+3}( x_{2}) +\frac{q^{2n+4}}{z^{2}}( x_{2}) -q^{2n+3}(
x_{_{1}}) -\frac{q^{2n+4}}{z^{2}}( x_{_{1}})
\\
&=&q^{2}(x_{2})\Big(q^{2n+1}(x_{2})+\frac{q^{2n+2}(x_{2})}{z^{2}}\Big)
-q^{2}(x_{_{1}})\Big(q^{2n+1}(x_{_{1}})+\frac{q^{2n+2}(x_{_{1}})}{z^{2}}
\Big) \\
&\geq &\Big(q^{2}(x_{2})-q^{2}(x_{_{1}})\Big)\Big(q^{2n+1}(x_{_{1}})+%
\frac{q^{2n+2}(x_{_{1}})}{z^{2}}\Big) \\
&\geq
&q^{2n+1}(x_{_{1}})\Big(q^{2}(x_{2})-q^{2}(x_{_{1}})\Big)\Big(1+
\frac{q(0)}{z^{2}}\Big) \geq 0.
\end{eqnarray*}%
Therefore, the function $x\mapsto
q^{2n+1}(x)+\frac{q^{2n+2}(x)}{z^{2}}$ is increasing in $[
0,x_{0}]$.
\end{proof}
We are now ready to prove Theorem \ref{the1}.
\begin{proof}  the help of the demonstration is to develop
$\dot{\theta}( x_{0},z)$ in entire series.\\ If $\varphi ( x_{0},z)
<\frac{\pi }{2},$ then the statement of the theorem immediately
follows from $\eqref{2.9} .$ If $\varphi( x_{0},z) \geq \frac{\pi
}{2},$ let $\varphi ( x_{0},z) =k\pi +\frac{\pi }{2}+D,$ with\
$0\leq D\leq \pi .$ Then
\begin{eqnarray*}
&&\int_{0}^{x_{0}}r^{2}(x)\frac{q(x)}{z}\Big(\sin ^{2}\varphi
-\varphi
\sin \varphi \cos \varphi\Big) dx \\
&=&\int_{0}^{\varphi ^{-1}( k\pi +\frac{\pi }{2}+D) }r^{2}(x)
\frac{q(x)}{z}\Big(\sin ^{2}\varphi -\varphi \sin \varphi \cos
\varphi\Big) dx \\
&=&\int_{0}^{\varphi ^{-1}( \frac{\pi }{2}) }[\tau]
dx+\sum_{i=0}^{k-1}\int_{\varphi ^{-1}( i\pi +\frac{\pi }{2})
}^{\varphi ^{-1}( ( i+1) \pi +\frac{\pi }{2}) }[\tau]
dx+\int_{\varphi ^{-1}( k\pi +\frac{\pi }{2}) }^{\varphi ^{-1}( k\pi
+\frac{\pi }{2}+D) }[\tau] dx,
\end{eqnarray*}
where $[\tau] =r^{2}(x)\frac{q(x)}{z}(\sin^{2}\varphi -\varphi \sin
\varphi \cos \varphi) .$ Since $\varphi\in]0,\frac {\pi}{2}],$ then
the first integral in the right hand of the last equality is
negative. In view of Lemma \ref{lem4}, together with \eqref{2.7}, it
follows
\begin{eqnarray*}
&&\int_{\varphi^{-1}( i\pi +\frac{\pi }{2}) }^{\varphi
^{-1}((i+1)\pi +\frac{\pi }{2})} r^{2}(x)\frac{q(x)}{z}\Big(\sin
^{2}\varphi -\varphi \sin \varphi \cos \varphi\Big) dx \cr &\leq& -
(i+1) \pi \int_{\varphi^{-1}( i\pi +\frac{\pi }{2}) }^{\varphi
^{-1}( (i+1)\pi +\frac{\pi
}{2})}r^{2}(x)\frac{%
q(x)}{z}\sin \varphi \cos \varphi dx \cr&&=-( i+1)\frac{\pi}{2}[
r^{2}]_{\varphi^{-1}( i\pi +\frac{\pi }{2}) }^{\varphi ^{-1}(
(i+1)\pi +\frac{\pi }{2})}.
\end{eqnarray*}
Since the logarithmic function is strictly increasing, it is enough
to prove
\begin{eqnarray}\label{aa}
&&-[\log ( r^{2}) ] _{\varphi ^{-1}( i\pi +\frac{ \pi }{2})
}^{\varphi ^{-1}( ( i+1) \pi +\frac{\pi }{2} ) }=-2\int_{\varphi
^{-1}( i\pi +\frac{\pi }{2}) }^{\varphi ^{-1}( ( i+1) \pi +\frac{\pi
}{2}) }\frac{r^{\prime }}{r}dx \cr&&=-2\int_{\varphi ^{-1}( i\pi
+\frac{\pi }{2}) }^{\varphi ^{-1}( ( i+1) \pi +\frac{\pi }{2})
}\frac{q(x)}{z}\sin \varphi \cos \varphi dx\leq 0.
\end{eqnarray}
By the change of variables $t=\varphi(x),$ we have
\begin{eqnarray*}
&&-\int_{\varphi ^{-1}( i\pi +\frac{\pi }{2}) }^{\varphi ^{-1}( (
i+1) \pi +\frac{\pi }{2}) }\frac{q(x)}{z}\sin \varphi \cos \varphi
dx =-\int_{i\pi +\frac{\pi }{2}}^{( i+1) \pi +\frac{\pi }{2}}\Big(\frac{%
\frac{q(\varphi ^{-1}( t) )}{z^{2}}\sin(t)\cos(t)}{1-\frac{%
q(\varphi ^{-1}( t) )}{z^{2}}\sin ^{2}(t)}\Big)dt.
\end{eqnarray*}
For $z^{2}\geq -q( 0) ,$ we have $\vert \frac{q(\varphi ^{-1}( t)
)}{z^{2}}\sin ^{2}(t)\vert <1$, and hence, the function
$\frac{1}{1-\frac{q(\varphi ^{-1}( t) )}{z^{2}}\sin ^{2}(t)}$ is
developable in entire series. Thus,

\begin{eqnarray*}
&&-\int_{i\pi+\frac{\pi}{2}}^{(i+1)\pi+\frac{\pi}{2}}
\Big(\frac{\frac{q(\varphi^{-1}(t))\sin(t)\cos(t)}{z^{2}}}
{1-\frac{q(\varphi^{-1}(t))\sin^{2}(t)}{z^{2}}}\Big)dt\cr&&=-\int_{i\pi
+\frac{\pi }{2}}^{( i+1) \pi +\frac{\pi }{2}}\sum_{n\geq 0}
\Big(\frac{q(\varphi^{-1}(t))}{z^{2}}\Big)^{n+1}\sin^{2n}(t)\sin(t)\cos(t)dt\cr
&&=-\sum_{n\geq 0}\frac{1}{z^{4n+2}}\int_{i\pi +\frac{\pi }{2}}^{(
i+1) \pi +\frac{\pi }{2}}
q^{2n+1}(\varphi^{-1}(t))\sin^{4n}(t)\sin(t)\cos(t)dt\cr &&-
\sum_{n\geq 0}\frac{1}{z^{4n+4}}\int_{i\pi +\frac{\pi }{2}}^{( i+1)
\pi +\frac{\pi }{2}}
q^{2n+2}(\varphi^{-1}(t))\sin^{4n+2}(t)\sin(t)\cos(t)dt\cr&&=-\sum_{n\geq
0}\frac{1}{z^{4n+2}}\int_{i\pi +\frac{\pi }{2}}^{( i+1)
\pi+\frac{\pi }{2}}
{\textstyle\Big(q^{2n+1}(\varphi^{-1}(t))+\frac{q^{2n+2}(\varphi^{-1}(t))}{z^{2}}
\Big)}\sin^{4n+3}\cos(t)dt\cr &&-\sum_{n\geq
0}\frac{1}{z^{4n+2}}\int_{i\pi +\frac{\pi }{2}}^{( i+1)
\pi+\frac{\pi }{2}}
q^{2n+1}(\varphi^{-1}(t))\sin^{4n+1}(t)\cos^{3}(t)dt.
\end{eqnarray*}
According to Lemma \ref{lem5} and as $z^{2}\geq -2q( 0) ,$ we get
\begin{eqnarray}\label{A}
&&-\int_{i\pi +\frac{\pi }{2}}^{( i+1) \pi}
\Big(q^{2n+1}(\varphi^{-1}(t))+\frac{q^{2n+2}(\varphi^{-1}(t))}{z^{2}}\Big)\sin^{4n+3}(t)\cos(t)dt
\cr &&-\int_{i\pi +\frac{\pi }{2}}^{( i+1) \pi}
q^{2n+1}(\varphi^{-1}(t))\sin^{4n+1}(t)\cos^{3}(t)dt\cr &&\leq
-{\textstyle\Big(q^{2n+1}(\varphi^{-1}(( i+1)
\pi))+\frac{q^{2n+2}(\varphi^{-1}(( i+1)
\pi))}{z^{2}}\Big)}\int_{i\pi +\frac{\pi }{2}}^{( i+1)
\pi}\sin^{4n+3}(t)\cos(t)dt \cr
&&-q^{2n+1}(\varphi^{-1}((i+1)\pi)\int_{i\pi +\frac{\pi }{2}}^{(
i+1) \pi}\sin^{4n+1}(t)\cos^{3}(t)dt\cr &&\leq \frac{1}{4n+4}
\Big[q^{2n+1}(\varphi^{-1}(( i+1)
\pi))+\frac{q^{2n+2}(\varphi^{-1}(( i+1) \pi))}{z^{2}}\Big]\cr
&&+\frac{1}{(4n+4)(2n+1)} q^{2n+1}(\varphi ^{-1}((i+1) \pi)).
\end{eqnarray}
In a similar way, we get
\begin{eqnarray}\label{B}
&&-\int_{( i+1) \pi}^{( i+1)
\pi+\frac{\pi}{2}}\Big[q^{2n+1}(\varphi^{-1}(t))+\frac{q^{2n+2}(\varphi^{-1}(t))}{z^{2}}\Big]\sin^{4n+3}(t)\cos(t)dt
\cr &&-\int_{( i+1) \pi}^{(
i+1)\pi+\frac{\pi}{2}}q^{2n+1}(\varphi^{-1}(t))\sin^{4n+1}(t)\cos^{3}(t)dt\cr
&&\leq- \frac{1}{4n+4}
\Big[q^{2n+1}(\varphi^{-1}((i+1)\pi))+\frac{q^{2n+2}(\varphi^{-1}((
i+1)\pi))}{z^{2}}\Big]\cr &&-\frac{1}{(4n+4)(2n+1)} q^{2n+1}(\varphi
^{-1}((i+1) \pi)).
\end{eqnarray}
Therefore, from \eqref{A} and \eqref{B} we obtain \eqref{aa}.\\
It is easily seen that if $0\leq D\leq \frac{\pi }{2}$, then
\begin{eqnarray*}
\int_{\varphi ^{-1}(k\pi +\frac{\pi }{2})}^{\varphi ^{-1}(k\pi
+\frac{\pi }{2}+D)}r^{2}(x)\frac{q(x)}{z}\Big(\sin ^{2}\varphi
-\varphi \sin \varphi \cos \varphi\Big) dx \leq 0.
\end{eqnarray*}%
If $\frac{\pi }{2}\leq D\leq \pi $, then
\begin{eqnarray*}
&&-\int_{\varphi^{-1}((k+1)\pi)}^{\varphi^{-1}(k\pi+\frac{\pi}{2}+D)}
r^{2}(x)\frac{q(x)}{z} \sin \varphi\cos\varphi dx\cr
&&\leq-\int_{\varphi^{-1}((k+1)\pi)}^{\varphi^{-1}((k+1)\pi+\frac{\pi}{2})}
r^{2}(x)\frac{q(x)}{z} \sin \varphi\cos\varphi dx,
\end{eqnarray*}  and hence,
\begin{eqnarray} \label{ACB}
&&\int_{\varphi^{-1}(k\pi+\frac{\pi}{2})}^{\varphi^{-1}(k\pi+\frac{\pi}{2}+D)}
r^{2}(x)\frac{q(x)}{z}\Big(\sin^{2}\varphi-\varphi\sin
\varphi\cos\varphi\Big) dx\cr
&&\leq-(k+1)\pi\int_{\varphi^{-1}(k\pi+\frac{\pi}{2})}^{\varphi^{-1}(k\pi+\frac{\pi}{2}+D)}
r^{2}(x)\frac{q(x)}{z} \sin \varphi\cos\varphi dx\cr&&
\leq-(k+1)\pi\int_{\varphi^{-1}(k\pi+\frac{\pi}{2})}
^{\varphi^{-1}((k+1)\pi+\frac{\pi}{2})} r^{2}(x)\frac{q(x)}{z} \sin
\varphi\cos\varphi dx \cr&&\leq0.
\end{eqnarray}
Therefore, by $\eqref{aa} $ and $\eqref{ACB},$ $ \dot{\theta}(
x_{0},z) \leq 0 $ for $z\geq \sqrt{-2q (0)}.$ Obviously,
 if there is a \\ $z\geq\sqrt{-2q(0)} $ with $\dot{\theta}(x_{0},z)
=0,$ then by Lemma \ref{lem4}, $q\equiv0$ in $[0,x_{0}]$. This
completes the proof of the theorem.
\end{proof}
\section{Proof Of Theorem \ref{the2}} \label{b}

Recall that $q(x)$ is monotone decreasing in $[x_{0},1].$ Let
$\tilde{q}(x)$ denotes the reverse of the potential, i.e.,
$\tilde{q}(x)=q(1-x)$. Then $\tilde{y}(x,z) $ is the solution of
System \eqref{H.S}, where $q(x)$ is replaced by $\tilde{q}(x)$. The
associated modified Pr\"{u}fer substitution is
\begin{eqnarray} \label{n3.2}
\left\{
\begin{array}{ll} \tilde{y}( x,z) =\tilde{r}( x,z) \sin (
z\tilde{\theta}( x,z) ),\\
\tilde{y}^{\prime }( x,z) =z\tilde{r}( x,z) \cos ( z\tilde{\theta}(
x,z) ),\\
\tilde{\theta}(0,z) =0.
 \end{array}
 \right.
\end{eqnarray}
As in \cite{5}, we have the following relations:
\begin{eqnarray} \label{3.2}
\left\{
 \begin{array}{ll} \tilde{y}( x,z_{n}) =( -1) ^{n+1}\frac{y(1-x,z_{n})}{
z_{n}r( 1,z_{n}) },\\
\tilde{r}( x,z_{n}) =\frac{r(1-x,z_{n})}{r( 1,z_{n}) },\\
\tilde{\theta}( x,z_{n}) =\frac{n\pi }{z_{n}}-\theta (1-x,z_{n}),
 \end{array}
 \right.
\end{eqnarray}
where $\lambda_{n}=z_{n}^{2}$ is an eigenvalue of Problem
\eqref{1.1}-\eqref{1.2}.

\begin{proof}{ of Theorem \ref{the2}}
As $\tilde{q}(x)=q(1-x),$ then $\tilde{q}(x)$ is monotone increasing
in $[0,1-x_{0}]$ and monotone decreasing in $[1-x_{0},1].$ Thus by
Theorem $\ref{the1},$ $ \tilde{\theta}( 1-x_{0},z) $ is decreasing
for $z\geq \sqrt{-2\tilde{q}(0)}=\sqrt{-2q( 1)}.$ Consequently, the
function $\Psi ( z) =\theta ( x_{0},z) +\tilde{\theta}( 1-x_{0},z) $
is decreasing for $z\geq \sqrt{-2q^{\ast }},$ where $q^{\ast }=\min
\{ q( 0) ,q( 1) \}$. Let $m$ be an integer such that $m<n$ and
$\lambda_{m}\geq-2q^{\ast}$. Then
\begin{eqnarray*}
\Psi ( \rho _{n}) =\frac{n\pi }{\rho _{n}}\leq \Psi ( \rho _{m})
=\frac{m\pi }{\rho _{m}},
\end{eqnarray*}
which implies that $ \frac{\lambda _{n}}{\lambda _{m}}\geq
\frac{n^{2}}{m^{2}}. $ On the other hand, if $|q^{\ast }|\leq
\frac{\pi^{2}}{3}$, then $$z_{1}=\sqrt{\lambda_{1}}\geq
\sqrt{\pi^{2}+q^{\ast }} \geq \sqrt{-2 q^{\ast }}.$$ Thus, in this
case $\lambda_{1}>0$ and $ \frac{\lambda _{n}}{\lambda _{m}}\geq
\frac{n^{2}}{m^{2}}$ for all $n>m\geq1.$ If equality holds, then
$\Psi ( z_{n}) =\Psi ( z_{m}) , $ so that $\dot{\Psi}(z) =0$ for
some $z>0.$ Thus $\dot{\theta}(x_{0},z)=\dot{\tilde{\theta}}(
1-x_{0},z) =0 ,$ and in view of Theorem \ref{the1}, $q\equiv0$ in $[
0,x_{0}] $ and $\tilde{q}\equiv0$ in $[ 0,1-x_{0}],$ i.e.,
$q\equiv0$ in $[0,1].$
\end{proof}

\begin{remark}\label{rem1}
~~~~~~~~~~~~~~~~~~~~~~~~~~~~~~~~~~~~~~~~~~~

                          It is easily seen that, if $q(x)$ is a nonnegative and
                             single-well potential, then $\dot{\theta}(x_{0},z)$ is developable as entire
                             series for $z>0$. Therefore,  the Pr\"{u}fer angle $\theta(x_{0},z)$
                             is monotone increasing for $z>0$. As a consequence,
                             $\frac{\lambda_{n}}{\lambda_{m}}\leq \frac{n^{2}}{m^{2}},$ for
                             $n>m\geq1.$

\end{remark}



\end{document}